\documentclass[reqno]{amsart}
\usepackage{amssymb}
\usepackage{graphics}
\usepackage{hyperref}

\def\Lie{\mathrm{Lie}}
\def\Perm{\mathrm{Perm}}

\def\ComTrias{\mathrm{ComTrias}}
\def\As{\mathrm{As}}
\def\PreLie{\mathrm{PreLie}}
\def\PostLie{\mathrm{PostLie}}
\def\idd{\mathrm{id}\,}
\def\Ker{\mathrm{Ker}\,}

\title{On the computation of Manin products for operads}

\author{V.~Yu.~Gubarev, P.~S.~Kolesnikov}

\address{Novosibirsk State University, Novosibirsk, Russia}
\email{vsevolodgu@mail.ru}

\address{Sobolev Institute of Mathematics, Novosibirsk, Russia}
\email{pavelsk@math.nsc.ru}

\begin{document}

\begin{abstract}
In the theory of binary quadratic operads, the white and black 
products of operads (called Manin products) play an important role. 
Given two such operads, the computation of either of their Manin 
products is a routine task.
We present and describe a computer program that helps to compute 
white and black Manin products of binary quadratic operads.
The same utility allows to find the Koszul-dual operads.
In particular, we compute 
the white product of the operads Lie and As (governing 
the varieties of Lie and associative algebras, respectively).
It turns out that the resulting operad is magmatic, i.e., 
defines the variety of all algebras with one bilinear operation.
\end{abstract}


\maketitle

\section{Introduction}
In a series recent works related to Leibniz algebras, Rota---Baxter operators, and dendriform 
algebras, the important role of the white and black products of operads (called Manin products \cite{GK1994}, 
they are denoted by $\circ $ and $\bullet$, respectively)
has been established. 
As a general reference in the operad theory one may apply to \cite{LV}.
However, most of examples that can be found in the literature
deal with white products $\mathcal P\circ \Perm$, $\mathcal P\circ \ComTrias$, 
and black products $\mathcal P\bullet \PreLie$, $\mathcal P\bullet \PostLie$. 
Here $\mathcal P$ is a binary quadratic operad, the definitions of 
$\Perm$, $\ComTrias $, $\PreLie $, and $\PostLie $ can be found in 
the reference preprint \cite{Zinbiel}.
Another class of operads has been considered in \cite{Stohmayer2008}, 
where white products also coincide with Hadamard products.

These examples of products are ``degenerate'' in the sense that 
for both $\Perm $ and $\ComTrias $  the white products
$\mathcal P\circ \Perm$, $\mathcal P\circ \ComTrias$ 
coincide with the Hadamard products
$\mathcal P\otimes \Perm$, $\mathcal P\otimes \ComTrias$
\cite{Vallette2008}. 
Since the white and black Manin products of binary quadratic operads 
are dual with respect to Koszul duality of operads, i.e., 
$(\mathcal P \bullet \mathcal Q)^! = \mathcal P^! \circ \mathcal Q^! $, 
and $\PreLie =\Perm^! $, $\PostLie = \ComTrias^! $, 
we have $\mathcal P\bullet \PreLie  = (\mathcal P^!\otimes \Perm)^!$,
$\mathcal P\bullet \PostLie  = (\mathcal P^!\otimes \ComTrias)^!$.


The computation of a Manin product may require 
a lot of computations even in the simplest cases. 
In a recent paper \cite{Loday11}, the following question 
has been stated: What is the black product 
of Com and As? (Equivalently, this is the dual to 
the white product of Lie and As.)

In this note, we describe in elementary terms the process 
of computation of Manin products for 
binary quadratic operads with no nontrivial unary operators.
We also describe an utility called {\tt manin} designed to perform 
the corresponding computations over the field of 
rational numbers. The source code of the utility is available at 
\url{http://math.nsc.ru/LBRT/a1/pavelsk/manin2.zip}.
It is designed by GNU Pascal version 20070904 based on gcc-4.1.3 20080704
for Ubuntu 2.1-4.1.2-29ubuntu2 (Free Software Foundation, Inc.). 
The executable file for Windows XP is compiled by Borland Delphi 5.0.  

\section{Computing Manin products}
Let $\Bbbk $ be a base field. 
If $\mathcal P$ is a variety of binary algebras over $\Bbbk $
defined by multilinear identities then $\mathcal P$
is governed by an operad which is also denoted by $\mathcal P$.
A binary quadratic operad corresponds to a variety 
whose defining identities have degrees 2 or 3. To define such an operad, 
we need a space of (binary) operations $E$ and the space of 
relations $R$, see \cite{GK1994} for details.

The space of multilinear terms of degree 3 is identified with 
\[
 E(3) = \Bbbk S_3\otimes _{\Bbbk S_2}(E\otimes E)
\]
where the action of $(12)\in S_2 $ on $E\otimes E$ is given by 
$\idd\otimes (12)$.

If $\mu $ is an
element of an $S_2$-module $E$, 
representing binary operation 
$(x_1,x_2)\mapsto x_1 * x_2$
then
$\mu^{(12)}$ corresponds to $(x_1,x_2)\mapsto x_2* x_1$.

If $\mu_1 $ and $\mu_2 $ represent two binary operations 
$(x_1,x_2)\mapsto x_1 *_k x_2$ ($k=1,2$)
then $\mu_1\otimes \mu_2 \in E\otimes E$
corresponds to the following function:
$(x_1,x_2,x_3)\mapsto (x_1 *_2 x_2)*_1 x_3 $.
It is convenient to identify 
the basic elements of $E(3)$ with appropriate monomials in
formal variables $x_1,x_2,x_3$ as stated in the table below.

\begin{center}
 \begin{tabular}{|c|c|}
\hline
$1\otimes_{\Bbbk S_2}(\mu_1\otimes \mu_2)$ &                                          $(x_1*_2 x_2)*_1 x_3$ \\
$1\otimes_{\Bbbk S_2}\big (\mu_1\otimes \mu_2^{(12)}\big )$ &                  $(x_2*_2 x_1)*_1 x_3$ \\
$1\otimes_{\Bbbk S_2}\big (\mu_1^{(12)}\otimes \mu_2\big )$ &                   $x_3*_1( x_1 *_2 x_2)$ \\
$1\otimes_{\Bbbk S_2}\big (\mu_1^{(12)}\otimes \mu_2^{(12)}\big )$ &       $x_3*_1( x_2 *_2 x_1)$ \\
\hline
$(13)\otimes_{\Bbbk S_2}(\mu_1\otimes \mu_2)$ &                                     $(x_3*_2 x_2)*_1 x_1$ \\
$(13)\otimes_{\Bbbk S_2}\big (\mu_1\otimes \mu_2^{(12)}\big )$ &              $(x_2*_2 x_3)*_1 x_1$ \\
$(13)\otimes_{\Bbbk S_2}\big (\mu_1^{(12)}\otimes \mu_2\big )$ &              $x_1*_1( x_3 *_2 x_2)$ \\
$(13)\otimes_{\Bbbk S_2}\big (\mu_1^{(12)}\otimes \mu_2^{(12)}\big )$ &   $x_1*_1( x_2 *_2 x_3)$ \\
\hline
$(23)\otimes_{\Bbbk S_2}(\mu_1\otimes \mu_2)$ &                                      $(x_1*_2 x_3)*_1 x_2$ \\
$(23)\otimes_{\Bbbk S_2}\big (\mu_1\otimes \mu_2^{(12)}\big )$ &               $(x_3*_2 x_1)*_1 x_2$ \\
$(23)\otimes_{\Bbbk S_2}\big (\mu_1^{(12)}\otimes \mu_2\big )$ &               $x_2*_1( x_1 *_2 x_3)$ \\
$(23)\otimes_{\Bbbk S_2}\big (\mu_1^{(12)}\otimes \mu_2^{(12)}\big )$ &    $x_2*_1( x_3 *_2 x_1)$ \\
\hline
\end{tabular}
\end{center}

Suppose $\mathcal P_1$ and 
$\mathcal P_2$ are two binary 
quadratic operads, 
and  $\mathcal P_i(1)= 1$, $i=1,2$.
Then $\mathcal P_i=\mathcal P(E_i, R_i)$, 
where $E_i$ are the spaces of binary operations 
considered as $S_2$-modules, 
$R_i$ are the spaces of quadratic relations.
Assume $\dim E_i<\infty $. 

Recall that the Hadamard product $\mathcal P=\mathcal P_1\otimes \mathcal P_2$
is given by the rule 
$\mathcal P(n)= \mathcal P_1(n)\otimes \mathcal P_2(n)$, $n\ge 1$, 
and the composition maps are expanded on $\mathcal P(n)$ 
in the componentwise way. It the same way, the structure of an 
$S_n$-module is defined on $\mathcal P(n)$: 
A permutation $\sigma \in S_n$ acts on $\mathcal P_1(n)\otimes \mathcal P_2(n)$
as $\sigma \otimes \sigma $.

By definition \cite{GK1994}, 
the white product $\mathcal P_1\circ \mathcal P_2$ is the sub-operad 
of $\mathcal P_1\otimes \mathcal P_2$ generated by the 
space of operations 
$E= E_1\otimes E_2$. It is known to be a binary quadratic operad. 
To compute the space of relations, 
consider the map 
\[
\tau: E(3) \to E_1(3)\otimes E_2(3)
\]
 defined by
\[
\tau: \sigma\otimes _{\Bbbk S_2}  ((a_1\otimes b_1)\otimes (a_2\otimes b_2))
\mapsto 
(\sigma\otimes _{\Bbbk S_2}(a_1\otimes a_2) )\otimes (\sigma\otimes _{\Bbbk S_2}(b_1\otimes b_2) ),
\]
$a_i\in E_1$, $b_i\in E_2$, $\sigma \in S_3$.
This is a well-defined $S_3$-linear map. Obviously, $\tau $ is injective. Denote by $D(E_1,E_2)$ the image of $\tau $.

Since $\mathcal P_i(3)=E_i(3)/R_i$ for $i=1,2$,
the desired space of relations (a subspace in $E(3)$) is exactly the kernel 
of 
\[
E(3) \overset{\tau}{\to}  E_1(3)\otimes E_2(3) \overset{\tau_1\otimes \tau_2}{\to} \mathcal P_1(3)\otimes \mathcal P_2(3), 
\]
where $\tau_i: E_i(3)\to \mathcal P_i(3)$ are the natural epimorphisms.
The kernel of $\tau_1\otimes \tau_2$ is equal to 
$R_1\otimes E_2(3) + E_1(3)\otimes R_2$. It remains to find the intersection 
of $\Ker (\tau_1\otimes \tau_2)$ with $D(E_1,E_2)$ and apply $\tau^{-1}$
to get $R$---the space of relations, defining $\mathcal P_1\circ \mathcal P_2$.

For every finite-dimensional $S_n$-module $M$, 
let $M^\vee $ stand for the dual space of $M$ considered as an
$S_n$-module 
with respect to sgn-twisted action: 
$\langle f^{\sigma }, e\rangle = -\langle f,e^{\sigma }\rangle $, 
$f\in M^\vee $, $e\in E$, $\sigma \in S_n$. 

Recall that if $\mathcal P = \mathcal P(E,R)$, $R\subseteq E(3)$, 
then the Koszul dual operad $\mathcal P^!$ is defined as $\mathcal P(E^\vee, R^\perp)$, 
where $E^\vee $ is the dual space to $E$ endowed with 
sgn-twisted $S_2$-action and $R^\perp $ is the subspace of $E^\vee (3)\simeq E(3)^\vee $
orthogonal to $R$. 

To get the black product of two binary quadratic operads $\mathcal P_1$ and $\mathcal P_2$, 
it is enough to compute 
\[
\mathcal P_1\bullet \mathcal P_2 = \big (\mathcal P_1^!\circ \mathcal P_2^!\big )^!.
\]

\section{Description of the program}

\subsection{Presenting initial data}
Each input operad $\mathcal P(E,R)$ is described in a separate file (e.g., the operad of Leibniz algebras---in the file {\tt leib}).
The description consists of several parts. The first line of the file contains the number $n=\dim E$.
The next $n$ lines present the action of $(12)\in S_2$ on $E$: Each $i$th line $x_{i1}\ x_{i2}\ \dots \ x_{in}$
consists of coordinates of $e_i^{(12)} = x_{i1}e_1 + x_{i2}e_2 + \dots + x_{in}e_n$. 
Below, the number of defining relation should be stated in a separate line 
 (followed by a comment, e.g., what an operad is defined by this file).
After that, the list of relations $R$ comes. They are presented as integer coordinate rows in the 
 following standard basis of $E(3)$:
\[
\begin{gathered}
a_1 = \idd\otimes (e_1\otimes e_1), \ \dots,\ a_n =  \idd\otimes (e_1\otimes e_n), \\
a_{n+1}=\idd\otimes (e_2\otimes e_1), \  \dots,\ a_{2n} =  \idd\otimes (e_2\otimes e_n), \\
\dotfill \\
a_{n(n-1)+1}=\idd\otimes (e_n\otimes e_1), \dots, \  a_{n^2} = \idd\otimes (e_n\otimes e_n),\\
a_{n^2+1}=(13)\otimes (e_1\otimes e_1), \  \dots,\  a_{2n^2}= (13)\otimes (e_n\otimes e_n), \\
a_{2n^2+1}=(23)\otimes (e_1\otimes e_1), \  \dots,\  a_{3n^2}= (23)\otimes (e_n\otimes e_n). \\
\end{gathered}
\]
Here we identify the space $E(3) = \Bbbk S_3\otimes _{S_2} (E\otimes E)$ with 
a sum of three copies of $E\otimes E$, i.e.,  
$E(3)\simeq V_3\otimes E\otimes E$, where $V_3$ is formally spanned by 
$\idd ,(13),(23)\in S_3$. 

The rest of the file can be used as a notebook, the program does not read these data.

For example, let us state the content of the file {\tt as} describing the operad governing 
associative algebras.

\medskip
\hrule
\smallskip

\verb! 2 % e_1 = x_1 x_2, e_2 = x_2 x_1!

\verb! 0 1 !

\verb! 1 0 !

\verb! 6 % ident. of associative algebra !

\verb! 1 0  0  0  0 0  0 -1 0 0  0 0 !

\verb! 0 0  0 -1  1 0  0  0 0 0  0 0 !

\verb! 0 0  0  0  0 0 -1  0 1 0  0 0 !

\verb! 0 0  0  0  0 1  0  0 0 0  0 -1 !

\verb! 0 0  1  0  0 0  0  0 0 -1 0 0 !

\verb! 0 -1 0  0  0 0  0  0 0 0  1 0 !

\verb! !

\verb! a_1=e \o e, a_2=e \o e^{12}, a_3=e^{(12)}\o e, a_4=e^{12}\o e^{12} !

\verb! a_5=(13)a_1, a_6=(13)a_2, a_7=(13)a_3, a_8=(13)a_4 !

\verb! a_9=(23)a_1, a_{10}=(23)a_2, a_{11}=(23)a_3, a_{12}=(23)a_4 !

\smallskip
\hrule
\medskip

\subsection{Usage of the utility}
To compute the white or black product of operads described in files {\tt file1} and {\tt file2}, 
type

\smallskip
{\tt manin w file1 file2} \quad or \quad {\tt manin b file1 file2}
\smallskip

\noindent
respectively. To compute the Koszul-dual operad $\mathcal P^!$ to an operad $\mathcal P$ described in {\tt file1}, type

\smallskip
{\tt manin d file1} 
\smallskip

The output is written into two files: {\tt result} and {\tt result.amx}. 
The first one contains a description of the resulting operad $\mathcal P(E,R)$ in the same  form 
as the input files do, i.e., after a minor editing (and, possibly, commenting) 
it can be used as an input file. The second one contains the description 
in a ``human-readable''  \AmS-TeX\ format. To write down the identities
one should assign binary operations to basic vectors of $E(3)$ and 
rewrite the relations in terms of these operations (as in the table stated above).
Note that the final form of identities highly depends on this assignment.

\subsection{Overview of the units}
The main program {\tt manin.pas} uses three units: {\tt lspace}, {\tt dynarr}, and {\tt shmidt}.

The unit {\tt lspace} contains the definitions of main types of data: {\tt Vector} and {\tt Space}. 
Vectors are presented as lists of integers, spaces---as lists of vectors. Also, in the unit
{\tt lspace} the main arithmetic operations with vectors and spaces are defined
(sum, tensor product, intersection), as well as input-output routines.

For intersection, the Hauss reduction method is implemented: If
$V_1$ is a $\mathbb Q$-linear span of vectors $a_1,\dots, a_m\in \mathbb Z^n$, 
$V_2$---of $b_1,\dots , b_l\in \mathbb Z^n$ then
the basis of $V_1\cap V_2$  can be found as follows. 
Consider a matrix of size $(k+l)\times (2n)$ given by
\[
\begin{pmatrix}
 \hbox{--- $a_1$ ---} & \hbox{--- $a_1$ ---} \\
 \hdots & \hdots \\
 \hbox{--- $a_m$ ---} & \hbox{--- $a_m$ ---} \\
 \hbox{--- $b_1$ ---} & \hbox{--- $0$ ---} \\
 \hdots & \hdots \\
\hbox{--- $b_l$ ---} & \hbox{--- $0$ ---} 
\end{pmatrix}
\]
and apply elementary transformations of rows to make it a trapezoid. 
All vectors remaining in the right half of the table, opposite to 
zero vectors in the left half, span the intersection $V_1\cap V_2$.

The unit {\tt dynarr} is just a description of two-dimensional dynamic arrays and 
procedures allowing to convert a list of vectors to an array and converse.

In the unit {\tt shmidt}, the orthogonalization procedure (Gram---Schmidt process) is implemented.
We use this process to compute the orthogonal complements.
Namely, in order to find the orthogonal component of a vector 
$v\in \mathbb Z^n$ relative to a subspace $V$ spanned by 
$a_1,\dots, a_m \in \mathbb Z^n$, 
we first make the vectors $a_1,\dots, a_m $ pairwise mutually orthogonal 
and then compute 
$v_1 = \langle a_1, a_1\rangle v - \langle a_1, v\rangle a_1$,
$v_2 = \langle a_2, a_2\rangle v_1 - \langle a_2, v_1\rangle a_2$,
and so on (cancellation of coefficients is applied on each step).
Thus, to find an orthogonal complement $V^\perp $
for $V\subset \mathbb Q^n$, we find  
orthogonal components of all standard basic vectors $e_1,\dots, e_n$
relative to $V$ and then apply Hauss reduction process. 

The Hauss reduction method as well as the Gram---Schmidt process are 
encoded in such a way that both input and output lists of vectors 
have integer coordinates.

\section{Examples}

\subsection{White product of Lie and Perm}
The white product of the operads governing Lie and Perm
 algebras is known to be the operad of Leibniz algebras \cite{Chap01}.
Entering

\smallskip
{\tt manin w lie perm}
\smallskip

\noindent
we obtain the following relations (written in {\tt result.amx}; here we have just replaced 
the \AmS-TeX\ commands
\verb!\pmatrix!
 and 
\verb!\endpmatrix!
with the corresponding \LaTeX\ environment):

\medskip
\hrule 
\smallskip

Space of operations $E$: $a_{1}, \dots , a_{2}$

$S_2$ acts by:
$$\begin{pmatrix}
0 & -1 \\
-1 & 0 \\
\end{pmatrix} $$
Relations:

$+1(23)\otimes_{S_2}(a_{1}\otimes a_{1}) -1(23)\otimes_{S_2}(a_{1}\otimes a_{2}) $

$+1(13)\otimes_{S_2}(a_{1}\otimes a_{1}) -1(13)\otimes_{S_2}(a_{1}\otimes a_{2}) $

$-1(id)\otimes_{S_2}(a_{1}\otimes a_{1}) +1(13)\otimes_{S_2}(a_{2}\otimes a_{2}) +1(23)\otimes_{S_2}(a_{2}\otimes a_{1}) $

$+1(id)\otimes_{S_2}(a_{2}\otimes a_{1}) -1(13)\otimes_{S_2}(a_{2}\otimes a_{1}) -1(23)\otimes_{S_2}(a_{1}\otimes a_{2}) $

$+1(id)\otimes_{S_2}(a_{1}\otimes a_{2}) -1(13)\otimes_{S_2}(a_{2}\otimes a_{2}) -1(23)\otimes_{S_2}(a_{2}\otimes a_{1}) $

$+1(id)\otimes_{S_2}(a_{2}\otimes a_{2}) -1(13)\otimes_{S_2}(a_{1}\otimes a_{2}) -1(23)\otimes_{S_2}(a_{2}\otimes a_{2}) $

\smallskip
\hrule
\medskip

Now, let us interpret $a_1$ as the operation $[x_1x_2]$. Then, according to the obtained $S_2$-action on $E$, 
$a_2=-a_1^{(12)}$, i.e., $a_2$ should be interpreted as $-[x_2x_1]$. Thus the six relations above turn into
\[
\begin{gathered}[]
[[x_1x_3]x_2]+[[x_3x_1]x_2], \\
[[x_3x_2]x_1]+[[x_2x_3]x_1], \\ 
-[[x_1x_2]x_3] + [x_1[x_2x_3]] - [x_2[x_1x_3]], \\
-[x_3[x_1x_2]] + [x_1[x_3x_2]] + [[x_3 x_1]x_2], \\
-[[x_2x_1]x_3] - [x_1[x_2x_3]] + [x_2[x_1x_3]], \\
[x_3[x_2x_1]] + [[x_2x_3]x_1] - [x_2[ x_3x_1]]. 
\end{gathered}
\]
These are corollaries of the only identity $[x[yz]]-[y[xz]] -[[xy]z]$, 
the left Leibniz identity.

\subsection{The black product of PreLie and As}
The operad of dendriform algebras \cite{Loday2001} is 
known to be the black product of operads 
governing the varieties of pre-Lie and associative algebras
(see also \cite{GubKol2012}).
The command 

\smallskip
{\tt manin b prelie as}
\smallskip

\noindent
generates the following output:

\medskip
\hrule
\smallskip

Space of operations $E$: $a_{1}, \dots , a_{4}$

$S_2$ acts by:
$$\begin{pmatrix}
0 & 0 & 0 & -1 \\
0 & 0 & -1 & 0 \\
0 & -1 & 0 & 0 \\
-1 & 0 & 0 & 0 \\
\end{pmatrix} $$
Relations:

$-1(id)\otimes_{S_2}(a_{4}\otimes a_{4}) +1(13)\otimes_{S_2}(a_{1}\otimes a_{1}) -1(13)\otimes_{S_2}(a_{1}\otimes a_{3}) $

$+1(id)\otimes_{S_2}(a_{3}\otimes a_{3}) -1(13)\otimes_{S_2}(a_{2}\otimes a_{2}) +1(13)\otimes_{S_2}(a_{2}\otimes a_{4}) $

$-1(id)\otimes_{S_2}(a_{4}\otimes a_{2}) +1(13)\otimes_{S_2}(a_{3}\otimes a_{1}) $

$+1(id)\otimes_{S_2}(a_{2}\otimes a_{2}) -1(id)\otimes_{S_2}(a_{2}\otimes a_{4}) -1(13)\otimes_{S_2}(a_{3}\otimes a_{3}) $

$+1(id)\otimes_{S_2}(a_{3}\otimes a_{1}) -1(13)\otimes_{S_2}(a_{4}\otimes a_{2}) $

$+1(id)\otimes_{S_2}(a_{1}\otimes a_{1}) -1(id)\otimes_{S_2}(a_{1}\otimes a_{3}) -1(13)\otimes_{S_2}(a_{4}\otimes a_{4}) $

$+1(13)\otimes_{S_2}(a_{4}\otimes a_{1}) +1(23)\otimes_{S_2}(a_{1}\otimes a_{1}) -1(23)\otimes_{S_2}(a_{1}\otimes a_{3}) $

$-1(id)\otimes_{S_2}(a_{4}\otimes a_{1}) -1(23)\otimes_{S_2}(a_{1}\otimes a_{2}) +1(23)\otimes_{S_2}(a_{1}\otimes a_{4}) $

$+1(id)\otimes_{S_2}(a_{3}\otimes a_{2}) +1(23)\otimes_{S_2}(a_{2}\otimes a_{1}) -1(23)\otimes_{S_2}(a_{2}\otimes a_{3}) $

$+1(13)\otimes_{S_2}(a_{3}\otimes a_{2}) +1(23)\otimes_{S_2}(a_{2}\otimes a_{2}) -1(23)\otimes_{S_2}(a_{2}\otimes a_{4}) $

$+1(13)\otimes_{S_2}(a_{4}\otimes a_{3}) +1(23)\otimes_{S_2}(a_{3}\otimes a_{1}) $

$-1(id)\otimes_{S_2}(a_{2}\otimes a_{1}) +1(id)\otimes_{S_2}(a_{2}\otimes a_{3}) -1(23)\otimes_{S_2}(a_{3}\otimes a_{2}) $

$-1(13)\otimes_{S_2}(a_{2}\otimes a_{1}) +1(13)\otimes_{S_2}(a_{2}\otimes a_{3}) +1(23)\otimes_{S_2}(a_{3}\otimes a_{3}) $

$-1(id)\otimes_{S_2}(a_{4}\otimes a_{3}) +1(23)\otimes_{S_2}(a_{3}\otimes a_{4}) $

$-1(id)\otimes_{S_2}(a_{1}\otimes a_{2}) +1(id)\otimes_{S_2}(a_{1}\otimes a_{4}) -1(23)\otimes_{S_2}(a_{4}\otimes a_{1}) $

$+1(13)\otimes_{S_2}(a_{3}\otimes a_{4}) +1(23)\otimes_{S_2}(a_{4}\otimes a_{2}) $

$+1(id)\otimes_{S_2}(a_{3}\otimes a_{4}) -1(23)\otimes_{S_2}(a_{4}\otimes a_{3}) $

$+1(13)\otimes_{S_2}(a_{1}\otimes a_{2}) -1(13)\otimes_{S_2}(a_{1}\otimes a_{4}) -1(23)\otimes_{S_2}(a_{4}\otimes a_{4}) $

\smallskip
\hrule
\medskip

The 18 relations above split into three orbits with respect to the action of $S_3$.
The representatives of these orbits are:

$-\idd \otimes_{S_2}(a_{4}\otimes a_{2}) +(13)\otimes_{S_2}(a_{3}\otimes a_{1}) $,

$-\idd \otimes_{S_2}(a_{4}\otimes a_{4}) +(13)\otimes_{S_2}(a_{1}\otimes a_{1}) - (13)\otimes_{S_2}(a_{1}\otimes a_{3}) $,

$\idd\otimes_{S_2}(a_{3}\otimes a_{3}) - (13)\otimes_{S_2}(a_{2}\otimes a_{2}) + (13)\otimes_{S_2}(a_{2}\otimes a_{4}) $.

Let us interpret $x_1\succ x_2$ as $a_1$ and $x_2\prec x_1$---as $a_2$.
Then  $a_3$ corresponds to $-x_1\prec x_2$ and $a_4$---to $-x_2\succ x_1$.
Hence, the defining identities of $\PreLie\bullet \As$ are:
\[
 \begin{gathered}
  x_3\succ (x_2\prec x_1) - (x_3\succ x_2)\prec x_1, \\
(x_3\succ x_2)\succ x_1 +(x_3\prec x_2)\succ x_1 - x_3\succ (x_2\succ x_1), \\
x_1\prec (x_2\prec x_3) + x_1\prec(x_2\succ x_3) - (x_1\prec x_2)\prec x_3.
 \end{gathered}
\]

\subsection{Black product of Com and As}
The command

\smallskip
{\tt manin b as comm}
\smallskip
\noindent

generates the following output:

\medskip
\hrule
\smallskip

Space of operations $E$: $a_{1}, \dots , a_{2}$

$S_2$ acts by:
$$\begin{pmatrix}
0 & -1 \\
-1 & 0 \\
\end{pmatrix} $$
Relations:

$+1(id)\otimes_{S_2}(a_{1}\otimes a_{1}) $

$+1(id)\otimes_{S_2}(a_{1}\otimes a_{2}) $

$+1(id)\otimes_{S_2}(a_{2}\otimes a_{1}) $

$+1(id)\otimes_{S_2}(a_{2}\otimes a_{2}) $

$-1(13)\otimes_{S_2}(a_{1}\otimes a_{1}) $

$-1(13)\otimes_{S_2}(a_{1}\otimes a_{2}) $

$-1(13)\otimes_{S_2}(a_{2}\otimes a_{1}) $

$-1(13)\otimes_{S_2}(a_{2}\otimes a_{2}) $

$-1(23)\otimes_{S_2}(a_{1}\otimes a_{1}) $

$-1(23)\otimes_{S_2}(a_{1}\otimes a_{2}) $

$-1(23)\otimes_{S_2}(a_{2}\otimes a_{1}) $

$-1(23)\otimes_{S_2}(a_{2}\otimes a_{2}) $

\smallskip
\hrule
\medskip

This is clear that such an operad corresponds to 
the variety of 3-nilpotent algebras. 
As a corollary (which is also easy to check by 
means of {\tt manin w lie as}), the white product 
$\Lie \circ \As $ is the magmatic algebra.

\end{document}